\documentclass{article}

\usepackage{amsfonts}
\usepackage{amsmath}
\usepackage{graphicx}
\usepackage{gensymb}
\usepackage{MnSymbol}
\usepackage{color}
\usepackage{lscape}
\usepackage{pdflscape}
\usepackage{pdfpages}
\usepackage{caption}
\usepackage{rotating}
\usepackage{tikz}
\usepackage{pgfplots}
\usepackage{wrapfig}
\usepackage{mathtools}
\usepackage{verbatim}
\usepackage{amsthm}
\usepackage{algorithm}
\usepackage{algpseudocode}

\newtheorem{proposition}{Proposition}[]

\newtheorem{conjecture}{Conjecture}[]

\DeclarePairedDelimiter{\ceil}{\lceil}{\rceil}
\usepackage[backend=bibtex,style=numeric,sorting=anyt, maxbibnames = 5]{biblatex}
\DeclareFieldFormat[article, inbook]{title}{#1}
\DeclareFieldFormat{pages}{#1}

\renewbibmacro{in:}{%
	\ifentrytype{article}{}{\printtext{\bibstring{in}\intitlepunct}}}
\title{\bf New Lower Bounds for van der Waerden Numbers \protect\\ Using Distributed Computing}
\author{\Large 
	Daniel Monroe \\
	7708 Hackamore Drive, Potomac MD 20854 \\
	Montgomery Blair High School \\
    {\tt JCMCC@danielmonroe.net}}

\begin{document}
	\maketitle
	
	\begin{abstract}
	This paper provides new lower bounds for van der Waerden numbers using Rabung's method, which colors based on the discrete logarithm modulo some prime. Through a distributed computing project with 500 volunteers over one year, we checked all primes up to 950 million, compared to 27 million in previous work. We point to evidence that the van der Waerden number for $r$ colors and progression length $k$ is roughly $r^k$.
	\end{abstract}
		\section{Introduction}

		Van der Waerden's Theorem \cite{Waerden} states that for any $k$ and $r$ there is an $n$ for which any $r$-coloring of $[n]$ admits a monochromatic arithmetic progression of length $k$, where $[n]$ is taken to denote $\{1, 2, \ldots, n \}$. The van der Waerden number $W(k, r)$ is the smallest such $n$.
		
		Only seven nontrivial van der Waerden numbers are known exactly. These are the entries in Table 1 without a \textgreater \space symbol. For other $k$ and $r$, only upper and lower bounds are available. 
		
		Van der Waerden's original proof gave bounds that were not primitive recursive. Shelah \cite{Shelah} gave a proof of bounds that were primitive recursive; in particular, he showed that $W(k,r)$ lies no higher than the fifth class of the Grzegorcyzk hierarchy.

		The current best upper bound is due to Gowers \cite{Gowers}, who proved that

		\begin{equation} \label{eq1} \notag
		\begin{split}
		W(k,r)<2^{\scriptscriptstyle {2^{r^{2^{2^{k+9}}}}}}.
		\end{split}
		\end{equation}
		
		There has also been work on upper bounds for $k = 3$ and $k = 4$. See Bloom and Sisack \cite{Bloom} for the latest result and a history in the former case and Green and Tao \cite{green_tao_2017} for the latter; neither paper explicitly gives bounds but bounds can easily be derived from them. See also \cite{Heule17}.
		
		There has also been some work on lower bounds. The construction of Berlekamp \cite{Berlekamp} shows that if $p$ is prime,
		
		 \begin{equation} \label{eq2} \notag
		\begin{split}
		W(p+1, 2)\ge p(2^p-1).
		\end{split}
		\end{equation}
		
		Lower bounds generated using the probabilistic method have also undergone significant development. The naive method gives a lower bound of $r^{k/2}$, which the Lov\'asz Local Lemma improves to $r^k/ek$. The best result of this type is due to Kuzik and Shabanov, who use a generalized Local Lemma with a recoloring algorithm in combination with a generalization of the Lov\'asz Local Lemma to establish the following:
		
		\begin{equation} \label{eq4} \notag
		\begin{split}
		W(k,r) \ge cr^k \text{ for some absolute constant } c.
		\end{split}
		\end{equation}

		Determining exact van der Waerden numbers seems to be computationally difficult. The last two discovered, $W(6,2)$ \cite{KourilP08} and $W(4,3)$ \cite{Kouril}, were found using SAT solvers running on special purpose computers. Those authors stated that $W(7,2)$ was not computable at this time, and perhaps never will be.
		
		Our knowledge of growth rates of van der Waerden numbers in general is extremely limited, though it is notable that all colorings known to be maximal take the form of the constructions of Rabung \cite{Rabung} and Herwig et al. \cite{HerwigHLM07}, which color based on the discrete logarithm modulo some prime, the latter being a structure-preserving transformation based on cyclic zippers. We provide numerical evidence on the growth rate of van der Waerden numbes based on such constructions and point to evidence that these are optimal.

		Through distributed computing, we applied 2 teraflops of computing power over 12 months, or about 500 CPU years on a 2GHz core, to Rabung's method. We checked the primes exhaustively up to 950 million, as compared to 27 million by Liang et al \cite{Xu12}. Those authors focused only on the 2-color case. Our methods are no different from those of previous work; the sole difference was the scale of the computation, which was distributed among the computers of 500 volunteers. We also used cyclic zippers, the method of the latter paper, which double the size of a coloring produced with Rabung's method, checking primes up to 40 million, without finding any new lower bounds superior to those of Rabung's method.\footnote{Our code and data are freely available at \url{https://github.com/hmonroe/vdw}.} It seems as though Rabung's method overtakes cyclic zippers for large $k$ as our computations would otherwise have shown a cyclic zipper generating a better lower bound for values of $W(k,2)$ with $k$ between 13 and 19 inclusive.
		
		The paper is organized as follows. Section 2 presents the constructions utilized. Section 3 gives evidence that the bounds generated are tight. Section 4 explains our computational infrastructure. Section 5 describes our results. Section 6 discusses the structure of the der Waerden number and presents three conjectures based on our findings.

		\section{Constructions}
		We used the constructions of discrete logarithms \cite{Rabung} and cyclic zippers \cite{HerwigHLM07}) with larger primes and substantially greater computing power; the algorithm for zipping can be found in the latter paper. We also include in Tables 1-3 bounds which follow from the recursion formulas of Xu \cite{Xu2013} and Blankenship et al. \cite{Blankenship}.
		
		Rabung's method colors $[p-1]$ so that the $n$th entry is assigned $\log_\rho n \mod r$, where $\rho$ is a fixed primitive root of $p$ and $\log_\rho n$ is the discrete logarithm. Notice that if $a, a+d, \ldots, a+(k-1)d$ are colored identically, then $ad^{-1}, ad^{-1}+1,\ldots, ad^{-1}+k-1$ with spacing 1 are as well. Therefore, modulo $p$, one only needs to check for progressions of spacing 1. Rabung's method can alternatively be seen as coloring with power residues.
		
		We $r$-color $0,1, \ldots, (k-1)p$ so that position $n$ is given color $\log_\rho n \mod r$ for $n$ not in the set ${0,p, \ldots, (k-1)p}$, then assign any colors to members of that set so that not all have the same color. Rabung showed that this coloring contains no monochromatic arithmetic progression of length $k$ if and only if: (a) there is no monochromatic arithmetic progression of spacing 1 in $1, \ldots, p-1$ and (b) if 1 and $p-1$ have the same color then $[(k-1)/2]$ do not all have the same color, while if 1 and $p-1$ have different colors then $[k-1]$ do not all have the same color \cite{Rabung}. Pseudocode is presented in Algorithm 1.

		The Cyclic Zipper Method  doubles the length of a coloring generated with the above method by interleaving it with itself. This was inspired by maximal colorings for $W(6,2)$, which were computed in \cite{KourilP08}.\footnote{In Table 2, lower bounds based on this method are marked Z. In some cases, \cite{HerwigHLM07} and \cite{RabungL12} were able to zip a  twice and quadruple its length. In Table 2, these are marked ZZ.} We checked $k$ up to 18 and $r$ up to 10 with cyclic zipping code shared by Rabung and Lotts \cite{RabungL12}. It was terminated after having checked all primes through 40 million, given CPU constraints, and found that improved bounds were no longer being produced. This method did recreate already known lower bounds for smaller $k$.
		
		Tables 1-3 include bounds found through the method of Xu \cite{Xu2013} of applying colorings recursively. He defines $WR(k, r)$, or the ring van der Waerden number, as the van der Waerden number over $\mathbb{Z}_p$ for some prime $p$. This is taken to mean 1 larger than the largest prime $p$ for which some $r$-coloring over $\mathbb{Z}_p$ contains no monochromatic $k$-length arithmetic progression. He shows that if $k\geq3, s\geq2, t\geq1, 5\leq n<WR(k, s)$, then $W(k, st)>p(W(k,t)-1)+1$. The case $t=1$ is particularly important; it lets one concatenate $k-1$ copies of a coloring of $\mathbb{Z}_p$. This is the concatenation which the power residue method uses, though it is able to add one color to the end. Note that throughout the rest of the paper the ring van der Waerden number will be taken to be that over $\mathbb{F}_p$ rather than $\mathbb{Z}_p$.
		
		Finally, in Tables 1-3, we used the recursion formula of Blankenship et al \cite{Blankenship}, who showed that $W(k,r)>p\cdot(W(k,r-\ceil{\frac{r}{p}})-1)$ whenever $p$ is prime and less than or equal to $k$. This generalizes Berlekamp's \cite{Berlekamp} lower bound for $r=2$.
		
		\section{Evidence for Optimality of Constructions}
		
		Our bounds may shed light on the growth of the van der Waerden Numbers. Five of the seven known exact van der Waerden numbers (excluding $W(3,2)$ and $W(3,3)$) have tight lower bounds based on the power residue method and Cyclic Zipping. These exceptions are likely due to a unique behavior among the van der Waerden numbers for $k=3$; see Heule's paper \cite{Heule17} for more details.
		
		In particular, these two methods have been shown to produce all maximal colorings for all but the smallest known two-color van der Waerden numbers.
		
		There is another reason to expect the power residue method to give good lower bounds. If $\log_\rho n = \log_\rho (n + 1) = \log_\rho (n + 2) \mod r$, which is an arithmetic progression of spacing 1 and length 3, then for any $c$, $\log_\rho cn = \log_\rho (cn + c) = \log_\rho (cn + 2c) \mod r$, which is an arithmetic progression of spacing $c$. Therefore, these sequences are saturated arithmetic progressions of length $k-1$ of every possible spacing, which is a reason to think they are good lower bounds for $W(k,r)$.
		
		However, there are cases where SAT solvers have outperformed the Rabung's method plus Cyclic Zipping; for $W(5,3)$ a SAT solver \cite{Heuleetal} beat those two methods.

		\section{Computations}
		We used Berkeley Open Infrastructure For Network Computing (BOINC) to distribute the work among volunteers' computers. Two teraflops of computing power were utilized over 12 months, or around 500 CPU years on a 2 GHz core. To validate the results, two computers applied the Rabung's method to each prime. There were a total of 516 volunteers and 1760 computers in 53 countries. The program had Linux and Windows versions and was written in C++. Each computer was assigned as input a range of 250 integers (for instance, from 1,000 to 1,249), identified primes in this range, and applied the Rabuung's method to them. They then, for each $r$ recorded the longest progression $k-1$ in that coloring, and on that basis reported a lower bound for $W(k,r)$. The server received and amalgamated these reports to record the best lower bounds as shown in Table 1. There was no other pre- or post-processing by the server; the main computational work was done by the client. For primes up to 40 million, the volunteers' computers also applied Cyclic Zipping to the generated colorings.
		
		Both memory and CPU time were bottlenecks. For a given prime $p$, to apply Rabung's method, the program populates a one-dimensional array of length $p-1$ with the powers of a primitive root. As an optimization of memory and CPU time, the same array is used for two to seven colors. To check a prime of nearly one billion, the array would have nearly one billion entries stored in RAM. The CPU time for to populate the array for each prime was 3-4 minutes on average, depending on the speed of the volunteer's computer. It would be possible to store the array more compactly if the focus were on $r=2$, i.e., with one bit per array entry, through bit stuffing, although we did not pursue that given the CPU constraint. There is no apparent way to improve on CPU time by populating the array more efficiently.
		
		There are a number of reasons that a reader should have confidence in our results reported in Tables 1-3. First, two computers verified each computation in the BOINC project. Secondly, we reproduced all known lower bounds found using these methods. Thirdly, the reader can verify that the primes in Table 2 give the lower bounds in Table 1 using our source code, available at https://github.com/hmonroe/vdw with and without zipping. That site provides versions compiled for Windows, and sample input and output files. Verifying that the primes shown in Table 2 give the \textit{best possible} lower bounds would require similar resources as our project, that is, around 500 CPU years on a 2 GHz core.

		\section{Current Bounds}
		We computed lower bounds for a wider range of $k$ and $r$ than did previous work, as can be seen in Tables 1-3. These new lower bounds are shown in bold. Bounds not in bold are the best known of previous work: Rabung and Lotts \cite{RabungL12}, Herwig et al \cite{HerwigHLM07}, Kouril and Paul \cite{KourilP08}, Landman et al \cite{Landman}, Landman and Robertson \cite{landmanrobertson}, Rabung \cite{Rabung}, and Liang et al \cite{Xu12}. We checked the number of colors, $r$, up to 10, and progression length, $k$, up to 25, though we only present bounds for $r$ at most 4 since for more colors our bounds were beaten by the recursions of Xu and Blankenship.
		
		Our lower bounds on $W(k,2)$ grow roughly exponentially in $k$. Letting $WR'(k,r)$ and $W'(k,r)$ be this paper's lower bounds on the ring van der Waerden (that being van der Waerden numbers on $\mathbb{F}_p$ rather than $\mathbb{Z}_p$) and van der Waerden numbers, the ratio $W'(k+1,2)/W'(k,2)$ seems to oscillate between 2 and 2.7 when $k>14$, which is shown in Figure 1. The ratio $W'(k+1,3)/W'(k,3)$ seems to hover around 3 or 4. 
		There are several lower bound formulas with this growth rate. Berlekamp's \cite{Berlekamp} bound states that $W(p+1, 2)\ge p(2^p-1)$ for primes $p$. Blankenship et al \cite{Blankenship} show that $W(p+1,r)>p^{r-1}2^p$ for primes $p$, which is a generalization of this bound. Landman and Robertson \cite{landmanrobertson} generalized it in a different direction, showing that for primes $q$, $p\ge5$, $W(p+1,q)\ge p(q^p-1)+1$. 
		
		The strength of power residue colorings manifests itself in smaller van der Waerden numbers which are known exactly. As can be seen in Table 2, Rabung's method and cyclic zippers produce the largest possible colorings for 5 of the 7 currently known nontrivial van der Waerden numbers, that being Rabung's method for $W(4,2), W(4, 3), W(3, 4)$ and the cyclic zipper method for $W(5, 2)$ and $W(6, 2)$. In fact, for each of these, the respective method produces all maximal colorings.
		
		\section{The Structure of van der Waerden Numbers}
		
		The strength of Rabung's method suggests an interplay between the van der Waerden number and the ring van der Waerden number. Xu's method of concatenation can be applied to bound the former substantially above the latter. We present a simple argument that runs the other way.
		
		\begin{proposition}
			If every $r$-coloring of the Galois Field $\mathbb{F}_p$ admits monochromatic arithmetic progressions of length $2k(k-1)$, then $W(k,r) \leq p$.
		\end{proposition}
		
		\begin{proof}
			Consider the remainders of $0, d, 2d, \ldots, (k-1)d$ modulo $p$. There is an $1\leq r \leq k-1$ so that $rd$ leaves a remainder at most $p/k$, for the distances between consecutive terms sum to $p$ and there are exactly $k$ of them, so we may just take the difference of the two terms.
			
			Let the arithmetic progression in $\mathbb{Z}_p$ be $a, a + d, \ldots, (2k(k-1)-1)d$. The progression $a, a + rd, a + 2rd, \ldots, a + (2k-2)rd$ wraps around at most once, so it may be split into two progressions, one of which has at least $k$ terms.
		\end{proof}
		
		We conjecture that the van der Waerden number for large $k$ grows like $r^{\text{O}(k^2)}$, and that it may even grow as slowly as $r^k$ based on the following argument. If the following conjecture were shown, the structure of the van der Waerden number would be nearly completely determined, barring the case of large $r$:
		
		\begin{conjecture}
			For any prime of the form $rt+1$, a residue coloring minimizes progression length among all $r$-colorings of $\mathbb{F}_p$, or the Galois field over $p$ elements.
		\end{conjecture} 
	
		 It was shown in \cite{peralta} that if $p > Ck2^k$ for some absolute constant $C$, there are $k$ consecutive quadratic residues modulo $p$. Assuming the above conjecture, the two-color van der Waerden number would grow no more quickly than exponentially in the square of $k$. Looking at Figure 2 gives strong intuition in the regard of arithmetic progressions of residues of squares and cubes; note that the low value at $k=25$ for 2 colors suggests that the bound here can be improved. The dips for even $k$ suggests that arithmetic progressions of quadratic residues tend to be of even length.
		
		\begin{conjecture}
			Let $\pi(k, r)$ he greatest $p$ for which $p = rt+1$ and $p$ avoids arithmetic progressions of length $k$ among $r^\text{th}$ power residues. Then $\pi(k, r) = \Theta(kr^k)$.
		\end{conjecture}

		If colorings generated with quadratic residues are optimal for two colors for sufficiently large $k$, $W(2,k)$ would grow somewhat like $k^2 2^k$. We do not believe, however, that the method is optimal for larger numbers of colors. The four-color van der Waerden number, for example, would grow at least as quickly as $k^3 4^k$, as opposed to a $k^2 4^k$ bound generated naively. The recursion of Xu would allow one to combine factors of $k$ from the two-color van der Waerden number and the corresponding ring van der Waerden number.
		
		We conjecture that the van der Waerden number grows roughly as follows:
		
		\begin{conjecture}
			$\underset{k \to \infty}{\lim}\dfrac{W(k, r)}{W(k - 1, r)} = r$.
		\end{conjecture}
		
		\section{Acknowledgments}
		We are indebted to John Rabung for providing us with advice while planning the work. We appreciate the guidance of Jay Cummings, Bill Gasarch, Ben Green, Stanislaw Radzizowski, Aaron Robertson, Vladislav Taranchuk, Xiadong Xu, and Daniel Zhu. Any remaining errors are our own. We would also like to thank the over 500 volunteers who contributed their computing resources, without whom it would not have been possible to carry out a computation of this magnitude.

		\newpage
		\begin{algorithm}[]
			
			\caption{Rabung’s method}\label{rabung}
			\begin{algorithmic}[1]

				\State Initialize array $x_1, x_2, \ldots, x_{p-1}$ \Comment{discrete logarithms}
				\State Find primitive root $\rho$ of $p$.
				\State $a\gets1$ \Comment{logarithm value}
				\State $i \gets 0$ \Comment{power}
				\While{$a \neq 1$}
				\State $a\gets a*\rho \mod p$
				\State $i\gets i+1$
				\State $x_a \gets i$
				\EndWhile
				
				\ForAll{$2 \leq r \leq 10$} \Comment{number of colors}
				\ForAll {$1 \leq j \leq p-1$} \Comment{create coloring}
				\State $c_j \gets x_j \mod r$
				\EndFor
				\State $k \gets$ longest monochromatic run in $c_i$
				\State $k_0 \gets$ longest monochromatic run starting with $c_1$
				
				\If{not $(c_1 = c_{p-1}) \land (c_1 = \ldots = c_{(k-1)/2})$ and not $(c_1 \neq c_{p-1}) \land (c_1 = \ldots = c_{(k-1)/2})	$}
 				\State $W'(k,r) \gets \max(W'(k,r), (k-1)p + 1)$ \Comment{update}
				
				\EndIf
				\EndFor

			\end{algorithmic}
		\end{algorithm}

		\newpage
		
		\begin{tikzpicture}[baseline={(current bounding box.center)}]
		\begin{axis}[
		title={Figure 1. Growth Rate of Lower Bounds $W'(k,r)$},
		xlabel={Length of Monochromatic Arithmetic Progression $k$},
		ylabel={$\dfrac{W'(k, r)}{W'(k - 1, r)}$},ylabel style={rotate=-90},
		xmin=8, xmax=25,
		ymin=0, ymax=6,
		xtick={9,11,13,15,17,19,21,23,25},
		ytick={0,2,4,6},
		xtick pos=left,
		legend pos=north east,
		ymajorgrids=true,
		grid style=dashed,
		every axis plot/.append style={thick},
		]
		
		\addplot[]
		coordinates {
			(8,3.1)
			(9,3.6)
			(10,2.5)
			(11,1.9)
			(12,3.3)
			(13,2.6)
			(14,1.9)
			(15,2.7)
			(16,1.9)
			(17,2.8)
			(18,2.0)
			(19,2.8)
			(20,2.1)
			(21,2.7)
			(22,1.8)
			(23,2.4)
			(24,1.77)
			(25,2.095)

		};
		\addplot[dashed]
		coordinates {
			(8,4.9)
			(9,3.9)
			(10,4.5)
			(11,4.5)
			(12,4.3)
			(13,3.2)
			(14,2.7)
			(15,3.4)
			(16,4.1)
			
		};
		\legend{2 colors,3 colors}
		
		\end{axis}
		\end{tikzpicture}
		
		\vspace{1in}
	
		\begin{tikzpicture}[baseline={(current bounding box.center)}]
		\begin{axis}[
		title={Figure 2. Growth Rate of Lower Bounds $WR'(k,r)$},
		xlabel={Length of Monochromatic Arithmetic Progression $k$},
		ylabel={$\dfrac{WR'(k, r)}{r^k}$},ylabel style={rotate=-90},
		xmin=4, xmax=25,
		ymin=0, ymax=35,
		xtick={4,6,8,10,12,14,16,18,20,22,24},
		ytick={0,10,20,30,40},
		xtick pos=left,
		legend pos=north west,
		ymajorgrids=true,
		grid style=dashed,
		every axis plot/.append style={thick},
		]
		
		\addplot[]
		coordinates {
			(3,.375)
			(4,.688)
			(5,1.16)
			(6,2.17)
			(7,4.82)
			(8,4.18)
			(9,6.62)
			(10,11.2)
			(11,9.72)
			(12,9.38)
			(13,16.7)
			(14,14.6)
			(15,18.8)
			(16,16.7)
			(17,22.1)
			(18,20.4)
			(19,26.5)
			(20,26.4)
			(21,33.6)
			(22,29.3)
			(23,33.6)
			(24,28.4)
			(25,28.6)

		};
		\addplot[dashed]
		coordinates {
			(4,1.20)
			(5,0.992)
			(6,2.44)
			(7,3.34)
			(8,5.19)
			(9,5.92)
			(10,7.85)
			(11,10.5)
			(12,13.5)
			(13,13.1)
			(14,10.8)
			(15,11.2)
			(16,14.2)
			
		};
		\legend{2 colors,3 colors}
		
		\end{axis}
		\end{tikzpicture}

		\newpage

			\begin{table}[]
				\footnotesize
				\centering
				\caption{Lower Bounds for $W(k, r)$}
				\label{my-label}
				\def\arraystretch{2}
				\begin{tabular}{lrrrrrrrrr}
					& 2 colors  & 3 colors & 4 colors                                                                           \\\hline
					Length 3  & 9                            & 27    & 76                                                        \\
					
					Length 4  & 35                           & 293                          & \textgreater1,048             &                \\
					Length 5  & 178                          & \textgreater2,173            & \textgreater17,705                     \\
					Length 6  & 1,132                        & \textgreater11,191           & \textgreater157,348           & \\
					Length 7  & \textgreater3,703            & \textgreater48,811           & \textgreater2,284,751              \\
					Length 8  & \textgreater11,495           & \textgreater238,400          & \textgreater12,288,155            \\
					Length 9  & \textgreater41,265           & \textgreater932,745          & \textgreater139,847,085                           \\
					Length 10 & \textgreater103,474          & \textgreater4,173,724        & \textgreater1,189,640,578                                                 \\
					Length 11 & \textgreater193,941          & \textgreater18,603,731       & \textgreater3,464,368,083                                      \\
					Length 12 & \textgreater638,727          & \textgreater79,134,144       & \textgreater37,054,469,451        \\
					Length 13 & \textgreater1,642,309        & \textbf{\textgreater251,282,317}      & \textgreater224,764,767,431            \\
					Length 14 & \textgreater3,118,350        & \textbf{\textgreater669,256,082}      &  \textgreater748,007,969,550                                            \\
					Length 15 & \textgreater8,523,047        & \textbf{\textgreater2,250,960,279}                              \\
					Length 16 & \textgreater16,370,086       & \textbf{\textgreater9,186,001,216}                                           \\
					Length 17 & \textgreater46,397,777 & \textbf{\textgreater15,509,557,937}          \\
					Length 18 & \textgreater91,079,252                        \\
					Length 19 & \textgreater250,546,915            \\
					Length 20 & \textgreater526,317,462       \\
					Length 21 & \textbf{\textgreater1,409,670,741}       \\
					Length 22 & \textbf{\textgreater2,582,037,634}            \\
					Length 23 & \textbf{\textgreater6,206,141,987} \\
					Length 24 & \textbf{\textgreater10,980,093,212} \\
					Length 25 & \textbf{\textgreater23,003,662,489} \\
					\hline
				\end{tabular}

				\textgreater \space refers to lower bounds. Bounds in bold are new or have been improved by this paper. We include results only up to 4 colors as the recursions of Xu and Berlekamp et al. dominated for higher $r$.
				
			\end{table}

		\begin{table}[]
			\footnotesize
			\centering
			\caption{Method Used to Find Lower Bounds for $W(k, r)$}
			\def\arraystretch{2}
			\begin{tabular}{lrrrrrrrrr}
				
				& 2 colors  & 3 colors & 4 colors                                                                \\\hline
				Length 3  &      && 37      \\
				Length 4
				& 11                          &              97                & 349                             \\
				Length 5
				& 11ZZ                       &                              & 2,213Z                    &             &    \\
				Length 6
				& 113Z                       &                                    \\
				Length 7
				& 617                       &                              &  3,703x617                              \\
				Length 8
				& 821Z                       & 34,057                    &     11,495x1,069               \\
				Length 9
				&                      & 116,593                   & 41,265x3,389                    \\
				Length 10
				& 11,497                    & 463,747                   & 103,474x11,497                 \\
				Length 11
				& 9,697Z                     & 1,860,373                 & 193,941x17,863            \\
				Length 12
				& 29,033Z                    & 7,194,013                 & 638,727x58,013        \\
				Length 13
				& 136,859                   & \textbf{20,940,193}                & 1,642,309x136,859        \\
				Length 14
				& 239,873                   & \textbf{51,481,237}                &  3,118,350x239,873                                       \\
				Length 15
				& 608,789                   & \textbf{160,782,877}     \\
				Length 16
				& 1,091,339                 & \textbf{612,400,081}    \\
				Length 17
				& 2,899,861 & \textbf{969,347,371}   \\
				Length 18
				& 5,357,603    \\
				Length 19
				& 13,919,273   \\
				Length 20 
				& 27,700,919     \\
				Length 21
				& \textbf{70,483,537}   \\
				Length 22
				& \textbf{122,954,173}    \\
				Length 23
				& \textbf{282,097,363}   \\
				Length 24 & \textbf{477,395,357} \\
				Length 25 & \textbf{958,485,937} \\
				\hline
			\end{tabular}

			The numbers above are the primes we used to find the lower bounds. Results with the Cyclic Zipper Method \cite{HerwigHLM07} are shown above. Z=zipped once, ZZ=zipped twice. Entries with two numbers separated by an ``x"  are the numbers used in Xu's \cite{Xu2013} method. Entries with a prime and van der 
			Waerden number separated by a $\cdot$ were produced with Blankenship et al's \cite{Blankenship} recurrence relation. The lower bounds that are bold are new or have been improved by this paper. For primes close to 1 billion, there may be better results that we did not check (as in the case of $W(25,2)$).
		\end{table}

		\begin{table}[]
			
			\footnotesize
			\centering
			\caption{References for Lower Bounds for $W(k, r)$}
			\def\arraystretch{2}
			\begin{tabular}{lrrrrrrrrr}
				& 2 colors  & 3 colors & 4 colors                                                                               \\\hline
				Length 3  & \cite{Chvatal}                            &  \cite{Chvatal}                           & \cite{Beeler}                   \\

				Length 4  &  \cite{Chvatal}                           &\cite{Kouril}                         & \cite{Rabung}              \\
				
				Length 5  &   \cite{Stevens_Shantaram}                      & \cite{Heuleetal}            & \cite{HerwigHLM07}           \\
				
				Length 6  & \cite{KourilP08}                       & \cite{Heuleetal}           & \cite{Xu2013}           &   \\

				Length 7  & \cite{Rabung}            & \cite{Delft}          & \cite{Xu2013}           \\
				
				Length 8  & \cite{HerwigHLM07}         & \cite{HerwigHLM07}          & \cite{Xu2013}          \\
				
				Length 9  & \cite{HerwigHLM07}           & \cite{RabungL12}          & \cite{Xu2013}         \\
				Length 10 & \cite{Rabung}         & \cite{RabungL12}        & \cite{Xu2013}     \\
				Length 11 & \cite{RabungL12}          & \cite{RabungL12}       & \cite{Xu2013}      \\
				Length 12 & \cite{RabungL12}          & \cite{RabungL12}      & \cite{Xu2013}     \\
				Length 13 & \cite{Xu12}        & *      & \cite{Xu2013}    \\
				Length 14 & \cite{Xu12}        & *     & \cite{Xu2013} \\
				Length 15 & \cite{Xu12} & * \\
				Length 16 & \cite{Xu12}       & *    \\
				Length 17 & \cite{Xu12} & *     \\
				Length 18 & \cite{Xu12}       \\
				Length 19 & \cite{Xu12}      \\
				Length 20 & \cite{Xu12}      \\
				Length 21 & *    \\
				Length 22 & *    \\
				Length 23 & *    \\
				Length 24 & *    \\
				Length 25 & *    \\
				
				\hline
			\end{tabular}
			
			The asterisks (*) represent our new bounds.
		\end{table}

		\newpage
		
		\clearpage
		\printbibliography
		
\end{document}